\documentclass[11pt]{article}

\usepackage[letterpaper,margin=2.5cm]{geometry}
\usepackage[T1]{fontenc}
\usepackage{lmodern}
\usepackage{microtype}
\usepackage{amsmath,amssymb,amsfonts,amsthm}
\usepackage{mathtools}
\usepackage{graphicx}
\usepackage{booktabs,multirow}
\usepackage{algorithm,algpseudocode}
\usepackage{subcaption}
\usepackage{pifont}

\usepackage[numbers,sort&compress]{natbib}
\usepackage[colorlinks,linkcolor=red,anchorcolor=blue,citecolor=green]{hyperref}
\hypersetup{
  pdftitle={CI-PINN: Causal Integral Physics-Informed Neural Network for Solving Evolution Equations},
  pdfauthor={Xiaodong Feng, Ziyu Sun, Tao Tang, Xiaoliang Wan, Tao Zhou},
  pdflang={en-US},
  pdfsubject={Research article for Computer Methods in Applied Mechanics and Engineering},
  pdfkeywords={Physics-informed neural networks, Temporal causality, Evolution equations, Sparse collocation}
}

\newcommand{\bey}{\begin{eqnarray}}
\newcommand{\eey}{\end{eqnarray}}

\newcommand{\beq}{\begin{equation}}
\newcommand{\eeq}{\end{equation}}
\numberwithin{equation}{section}

\theoremstyle{plain}

\newcommand{\cmark}{\ding{51}}
\newcommand{\xmark}{\ding{55}}

\graphicspath{{figs/}}

\date{}
\title{CI-PINN: Causal Integral Physics-Informed Neural Network for Solving Evolution Equations}
\author{
Xiaodong Feng\thanks{Institute for Advanced Study, Beijing Normal-Hong Kong Baptist University,  Zhuhai, Guangdong, China. \texttt{xiaodongfeng@bnbu.edu.cn}}
\and
Ziyu Sun \thanks{School of Mathematical Sciences, Beijing Normal University, Beijing 100875, China. \texttt{zysun@mail.bnu.edu.cn}}
\and
Tao Tang\thanks{Institute for Advanced Study, Beijing Normal-Hong Kong Baptist University,  Zhuhai, Guangdong, China; and School of Mathematics and Statistics, Guangzhou Nanfang College, Guangzhou, Guangdong, China. \texttt{ttang@bnbu.edu.cn}}
\and
Xiaoliang Wan \thanks{Department of Mathematics and Center for Computation and Technology, Louisiana State University, Baton Rouge 70803, USA. \texttt{xlwan@lsu.edu}}
\and
Tao Zhou \thanks{SKLMS \& Institute of Computational Mathematics and Scientific/Engineering Computing, Academy of Mathematics and Systems Science, Chinese Academy of Sciences. \texttt{tzhou@lsec.cc.ac.cn}}
}

\begin{document}
\maketitle

\begin{abstract}
  Physics-informed neural networks (PINNs) solve partial differential equations (PDEs) by incorporating governing physical laws into the training loss. For evolution equations, however, their conventional pointwise space--time representation does not explicitly encode temporal dependence, which can hinder accurate prediction. To mitigate this limitation, this work proposes a novel neural architecture termed a causal integral neural network (CinNet). The core module of CinNet is a Volterra-type causal integral term, which aggregates historical features to encode temporal dependence, thereby incorporating temporal causality at the architectural level rather than through training-level modifications as in many existing methods. Building on CinNet, we further develop a causal integral physics-informed neural network (CI-PINN) for solving evolution equations. Extensive numerical experiments on benchmark evolution equations demonstrate that the presented method outperforms various baseline PINN variants in terms of solution accuracy, with pronounced superiority under sparse-collocation scenarios. Additional empirical analyses show that CI-PINN exhibits low sensitivity to hyperparameter choices, while ablation studies confirm the effectiveness of the proposed network components.
\end{abstract}

\noindent\textbf{Keywords:}
Scientific machine learning, physics-informed neural networks, temporal causality, causal integral neural network, sparse collocation

\section{Introduction}\label{sec:intro}

In recent years, deep learning methods, especially physics-informed neural networks (PINNs)~\cite{raissi2019Physicsinformed}, have shown great potential for solving partial differential equations (PDEs), including evolution equations. PINNs approximate the solution by a neural network and incorporate the governing equations, initial conditions, and boundary conditions into the training objective. This formulation has several appealing advantages, such as mesh-free representation, automatic differentiation of differential operators, flexible incorporation of observational data, the ability to evaluate the learned solution at arbitrary space--time points, and potential applicability to high-dimensional problems.

Despite these advantages, standard PINNs still face serious difficulties when solving time-dependent PDEs. In the vanilla formulation, the solution is usually parameterized as a pointwise mapping from space--time coordinates to the solution value,
\begin{equation*}
  (t,\mathbf{x}) \longmapsto \hat{u}(t,\mathbf{x};\theta).
\end{equation*}
In this representation, time is treated as an additional input coordinate, and the prediction at a later time is not structurally generated from earlier states. This is inconsistent with the intrinsic temporal causality of evolution equations, where the state at time $t$ is determined by the initial condition and the evolution history before $t$. As a result, standard PINNs may fail to propagate initial information effectively and may converge to low-residual but dynamically inaccurate solutions, especially in long-time prediction, stiff dynamics, or sparse-collocation regimes.

To improve the temporal behavior of neural solvers for time-dependent PDEs, existing studies have explored several directions. One class of methods reformulates the problem through domain decomposition~\cite{jagtap2020Extended,meng2020PPINN}, sequential training~\cite{wight2021Solving,krishnapriyan2021Characterizing,mattey2022Novel,penwarden2023Unified,roy2024Exact}, or discrete time-stepping formulations~\cite{jung2024CEENs}. Another class modifies the training objective or optimization procedure, for example through curriculum learning~\cite{krishnapriyan2021Characterizing}, pretraining~\cite{guo2023Pretraining}, causally biased sampling~\cite{daw2023Mitigating}, causal loss weighting~\cite{wang2024Respecting}, directional irreversibility regularization~\cite{chen2025Enforcing}, or integral-form residual regularization~\cite{feng2026Integral}. There are also problem-formulation methods that exploit special structures of the underlying PDE, such as characteristic or Lagrangian coordinates for convection-dominated problems~\cite{mojgani2023Kolmogorov}. More recently, representation-level methods have been proposed to improve temporal information propagation, including space--time separated representations~\cite{gu2022Deep,feng2024Hybrid}, recurrent~\cite{ren2022PhyCRNet}, Transformer-based~\cite{zhao2024PINNsFormer}, state-space~\cite{xu2025SubSequential}, autoregressive~\cite{nagda2025PIANO}, kernel-based~\cite{su2025SPIKE}, and parameter-evolution architectures~\cite{du2021Evolutional,bruna2024Neural,chen2024TENG}. These methods provide important ways of incorporating temporal structure into neural PDE solvers.

In this work, we focus on methods with continuous space--time neural representations, because they retain the main advantages of PINNs: the learned solution is a continuous function of $(t, \mathbf{x})$, and the required derivatives can be computed by automatic differentiation wherever they exist, scattered space--time data can be naturally incorporated, and the solution can be evaluated at arbitrary time points without a discrete autoregressive rollout. Within this setting, several representative approaches have been developed to alleviate the temporal training difficulty of PINNs. PPINN~\cite{meng2020PPINN} decomposes a long-time evolution problem into multiple short-time subproblems, where a fast coarse solver provides sequential predictions and independent fine PINNs correct the solution in parallel. Wight and Zhao~\cite{wight2021Solving} proposed adaptive sampling and time-marching strategies to reduce optimization difficulties and improve the accuracy of PINNs for phase-field equations. Krishnapriyan et al.~\cite{krishnapriyan2021Characterizing} introduced curriculum and sequence-to-sequence training strategies to mitigate failure modes of PINNs in time-dependent problems. Mattey and Ghosh~\cite{mattey2022Novel} proposed a backward-compatible PINN that solves time-dependent PDEs sequentially over successive time segments using a single neural network, while penalizing deviations from previously learned solutions to preserve consistency over earlier time intervals. Guo et al.~\cite{guo2023Pretraining} developed a pretraining strategy that first trains PINNs on short early-time intervals and then uses the learned parameters and pseudo-labels to initialize and regularize full-domain training for challenging evolution PDEs. Penwarden et al.~\cite{penwarden2023Unified} proposed a unified causal sweeping framework that combines temporal domain decomposition, window-based collocation propagation, and transfer initialization to enforce causal information flow from early to later times. Causal PINNs~\cite{wang2024Respecting} enforce temporal causality through loss reweighting, so that earlier-time residuals are reduced before later-time residuals are emphasized. IR-PINNs~\cite{feng2026Integral} introduce an additional integral-form residual term into the loss function, which acts as a regularization mechanism and enhances temporal correlation during training.

Although these methods significantly improve the performance of PINNs for evolution equations, most of them introduce temporal causality indirectly. Time-marching and domain-decomposition methods rely on splitting the time interval and solving a sequence of subproblems, while causal training, pretraining, and integral-residual methods mainly modify the training objective or optimization process. In these approaches, the underlying neural representation is still essentially the standard pointwise space--time mapping $\hat{u}(t,\mathbf{x};\theta)$. Therefore, temporal causality is imposed mostly through the training procedure rather than being embedded directly into the solution representation itself. This motivates us to design a neural architecture that retains the continuous space--time PINN formulation while introducing causal historical dependence at the representation level.

Motivated by this observation, we propose a new network suitable for evolution systems, known as a causal integral neural network (CinNet), which introduces a history-dependent integral term and a learnable gating mechanism into each network layer. The representation at time $t$ explicitly depends on historical features over $[t_0,t]$, where $t_0$ is the initial time, thereby introducing a causal inductive bias into the solution approximation. Building on this architecture, we develop a causal integral physics-informed neural network (CI-PINN) for solving time-dependent PDEs without requiring time-domain decomposition or discrete autoregressive rollout. Numerical experiments show consistent accuracy improvements over standard PINN and competitive performance relative to Causal PINN. The sparse-collocation results suggest an implicit regularization effect, while sensitivity and ablation studies assess the quadrature resolution, loss weighting, and contributions of the main architectural components.

The remainder of this paper is organized as follows. In Section~\ref{sec:preliminaries}, we introduce the standard formulation of PINNs and illustrate one of their most important failure modes in solving evolution equations through numerical experiments on the classical Allen--Cahn equation. In Section~\ref{sec:methodology}, we present the proposed CinNet architecture, discuss the roles of its components, provide the practical implementation of CinNet, and describe the CI-PINN framework for solving evolution equations. In Section~\ref{sec:results}, we conduct numerical experiments to evaluate the performance of CI-PINN on several representative evolution equations by comparing it with baseline methods and analyze the results. Finally, Section~\ref{sec:conclusion} concludes the paper and discusses future research directions.

\section{Preliminaries}
\label{sec:preliminaries}

This section first reviews the standard formulation of physics-informed neural networks for nonlinear evolution equations. We then discuss the inherent temporal causality of evolution problems and illustrate a representative failure mode of PINNs, which motivates the causal integral representation developed in the next section.

\subsection{Physics-informed neural networks}
\label{subsec:pinns}

We consider nonlinear time-dependent partial differential equations of the form
\begin{equation}\label{eq:pde}
  \begin{alignedat}{2}
     & u_t(t,\mathbf{x}) + \mathcal{N}[u](t,\mathbf{x}) = 0, &  & \qquad (t,\mathbf{x}) \in [t_0,T]\times\Omega,          \\
     & u(t_0,\mathbf{x}) = g(\mathbf{x}),                    &  & \qquad \mathbf{x} \in \Omega,                           \\
     & \mathcal{B}[u](t,\mathbf{x}) = h(t,\mathbf{x}),       &  & \qquad (t,\mathbf{x}) \in [t_0,T]\times \partial\Omega,
  \end{alignedat}
\end{equation}
where $u(t,\mathbf{x})$ is the unknown solution, $\Omega\subset\mathbb{R}^d$ is the spatial domain, and $[t_0,T]$ is the time interval. The operators $\mathcal{N}$ and $\mathcal{B}$ are the nonlinear spatial differential operator and the boundary operator, respectively, while $g$ and $h$ prescribe the initial and boundary conditions.

A multilayer perceptron (MLP) constructs a pointwise mapping from the space--time coordinates to the solution value through a sequence of affine transformations and nonlinear activation functions. An $L$-layer MLP can be written as
\begin{equation}\label{eq:mlp}
  \begin{aligned}
    \mathbf{y}^{\left( 0 \right)}\left( t,\mathbf{x} \right) & =\left( t,\mathbf{x} \right) ,                                                                                                             \\
    \mathbf{y}^{(\ell)}\left( t,\mathbf{x} \right)           & =\sigma \left( \mathbf{y}^{(\ell-1)}\left( t,\mathbf{x} \right) \mathbf{W}^{(\ell)}+\mathbf{b}^{(\ell)} \right) , \quad \ell=1,\ldots,L-1, \\
    \mathbf{y}^{\left( L \right)}\left( t,\mathbf{x} \right) & =\mathbf{y}^{\left( L-1 \right)}\left( t,\mathbf{x} \right) \mathbf{W}^{\left( L \right)}+\mathbf{b}^{\left( L \right)}.
  \end{aligned}
\end{equation}
Here, $\mathbf{y}^{(\ell)}(t,\mathbf{x})\in \mathbb{R} ^{D_{\ell}}$ is the output of the $\ell$-th layer, $\mathbf{W}^{(\ell)}\in \mathbb{R} ^{D_{\ell-1}\times D_\ell}$ and $\mathbf{b}^{(\ell)}\in \mathbb{R} ^{D_\ell}$ are the weights and biases of the $\ell$-th layer, and $\sigma(\cdot)$ is the activation function. The output of the MLP is given by $\hat{u}(t,\mathbf{x};\theta)=\mathbf{y}^{(L)}(t,\mathbf{x})$, where $\theta=\{\mathbf{W}^{(\ell)},\mathbf{b}^{(\ell)}: \ell=1,\ldots,L\}$ collects all trainable parameters.

A physics-informed neural network (PINN)~\cite{raissi2019Physicsinformed} approximates the solution of~\eqref{eq:pde} by the neural surrogate $\hat{u}(t,\mathbf{x};\theta)$. The network parameters are determined by minimizing an empirical loss consisting of the PDE residual loss and the penalties associated with the initial and boundary conditions:
\begin{equation}\label{eq:loss}
  \mathcal{L}(\theta) = \omega_{\mathrm{res}}\mathcal{L}_{\mathrm{res}}(\theta) + \omega_{\mathrm{ic}}\mathcal{L}_{\mathrm{ic}}(\theta) + \omega_{\mathrm{bc}}\mathcal{L}_{\mathrm{bc}}(\theta),
\end{equation}
where $\omega_{\mathrm{res}}$, $\omega_{\mathrm{ic}}$, and $\omega_{\mathrm{bc}}$ are nonnegative weights balancing the three loss components. The PDE residual associated with $\hat{u}$ is defined by
\begin{equation}\label{eq:residual}
  r(t,\mathbf{x};\theta) = \hat{u}_t(t,\mathbf{x};\theta) + \mathcal{N}[\hat{u}](t,\mathbf{x};\theta).
\end{equation}
Given the residual collocation points $\mathcal{D}_{\mathrm{res}}=\{(t_{\mathrm{res}}^i,\mathbf{x}_{\mathrm{res}}^i)\}_{i=1}^{N_{\mathrm{res}}}$,
initial points
$\{\mathbf{x}_{\mathrm{ic}}^i\}_{i=1}^{N_{\mathrm{ic}}}$,
and boundary points
$\{(t_{\mathrm{bc}}^i,\mathbf{x}_{\mathrm{bc}}^i)\}_{i=1}^{N_{\mathrm{bc}}}$,
the loss components are defined by
\begin{equation}\label{eq:loss_terms}
  \begin{aligned}
    \mathcal{L}_{\mathrm{res}}(\theta) & = \frac{1}{N_{\mathrm{res}}}\sum_{i=1}^{N_{\mathrm{res}}}\left|\hat{u}_t(t_{\mathrm{res}}^i,\mathbf{x}_{\mathrm{res}}^i;\theta)+\mathcal{N}[\hat{u}](t_{\mathrm{res}}^i,\mathbf{x}_{\mathrm{res}}^i;\theta)\right|^2, \\
    \mathcal{L}_{\mathrm{ic}}(\theta)  & = \frac{1}{N_{\mathrm{ic}}} \sum_{i=1}^{N_{\mathrm{ic}}} \left| \hat{u}(t_0,\mathbf{x}_{\mathrm{ic}}^i;\theta) - g(\mathbf{x}_{\mathrm{ic}}^i) \right|^2,                                                             \\
    \mathcal{L}_{\mathrm{bc}}(\theta)  & = \frac{1}{N_{\mathrm{bc}}} \sum_{i=1}^{N_{\mathrm{bc}}} \left| \mathcal{B}[\hat{u}] (t_{\mathrm{bc}}^i,\mathbf{x}_{\mathrm{bc}}^i;\theta) - h(t_{\mathrm{bc}}^i,\mathbf{x}_{\mathrm{bc}}^i) \right|^2.
  \end{aligned}
\end{equation}
In practice, the network parameters are determined by approximately minimizing $\mathcal{L}(\theta)$ using a gradient-based optimizer, and the resulting parameter vector $\theta^\star$ defines the trained PINN approximation $\hat{u}(t,\mathbf{x};\theta^\star)$.

The formulation above provides a flexible mesh-free framework for solving~\eqref{eq:pde}. Nevertheless, the MLP in~\eqref{eq:mlp} accepts time and space simply as coordinate inputs and does not explicitly encode the directional temporal structure of an evolution equation. The implications of this point are discussed next.

\subsection{Temporal causality and a failure mode of PINNs}
\label{sec:pinn_temporal_failure}

To make the temporal dependence explicit, we integrate~\eqref{eq:pde} over $[t_0,t]$ and use the initial condition to obtain
\begin{equation}\label{eq:evolution_integral_form}
  u(t,\mathbf{x})=g(\mathbf{x})-\int_{t_0}^{t} \mathcal{N}[u](s,\mathbf{x})\,\mathrm{d}s.
\end{equation}
Equation~\eqref{eq:evolution_integral_form} shows that the solution at time $t$ is determined by the initial condition together with the preceding evolution over $[t_0,t]$. Errors introduced near the initial time may therefore influence the entire subsequent trajectory.

By contrast, a conventional PINN represents the solution through the pointwise coordinate mapping $(t,\mathbf{x})\mapsto\hat{u}(t,\mathbf{x};\theta)$. Although predictions at different times are coupled through the shared network parameters, the representation at time \(t\) is not explicitly constructed from the preceding evolution. Consequently, satisfying the residual constraints at a finite set of space--time collocation points does not necessarily recover the correct temporal trajectory, particularly when temporal data are limited.
A standard energy estimate characterizes how the temporal distribution of the residual enters the error bound. Define
\begin{equation*}
  \mathcal{E}(t,\mathbf{x};\theta)=\hat{u}(t,\mathbf{x};\theta)-u(t,\mathbf{x}), \qquad \widetilde{\mathcal{L}}_{\mathrm{res}}(t;\theta)=\frac{1}{|\Omega|}\left\|r(t,\cdot;\theta)\right\|_{L^2(\Omega)}^2.
\end{equation*}
Under standard regularity assumptions and a one-sided $L^2(\Omega)$ stability condition on the evolution operator $-\mathcal{N}$, a Gr\"onwall argument gives the error estimate~\cite{zhang2026Solving}
\begin{equation}\label{eq:error_propagation}
  \begin{aligned}
    \left\|\mathcal{E}(t,\cdot;\theta)\right\|_{L^2(\Omega)}^2
    \leq{} & e^{(\alpha+1)(t-t_0)}\left\|\mathcal{E}(t_0,\cdot;\theta)\right\|_{L^2(\Omega)}^2 + |\Omega|\int_{t_0}^{t}e^{(\alpha+1)(t-s)}\widetilde{\mathcal{L}}_{\mathrm{res}}(s;\theta)\,\mathrm{d}s.
  \end{aligned}
\end{equation}
Here $\alpha\geq 0$ denotes the constant in the associated one-sided $L^2(\Omega)$ stability bound. Estimate~\eqref{eq:error_propagation} shows that the prediction error at time $t$ is controlled by the initial mismatch and a time-weighted accumulation of the residual over $[t_0,t]$. Since the contribution of the residual at time $s$ is weighted by $e^{(\alpha+1)(t-s)}$, its effect on the bound depends on when it occurs. Thus, the temporal distribution of the residual is relevant in addition to its global magnitude.

To examine this behavior numerically, we consider the Allen--Cahn equation~\eqref{eq:ac} and compute the following time-resolved diagnostics on the spatial test grid:
\begin{align*}
  \mathcal{L}_{\mathrm{res}}(t;\theta)
  = \frac{1}{N_x}\sum_{j=1}^{N_x}\left|r(t,x_j;\theta)\right|^2, \quad
  \mathrm{RL2E}(t)
  = \frac{\left(\sum_{j=1}^{N_x}\left|u(t,x_j)-\hat{u}(t,x_j;\theta)\right|^2\right)^{1/2}}{\left(\sum_{j=1}^{N_x}|u(t,x_j)|^2\right)^{1/2}}.
\end{align*}
\begin{figure}[!htbp]
  \centering
  \begin{subfigure}{0.45\textwidth}
    \includegraphics[width=\linewidth]{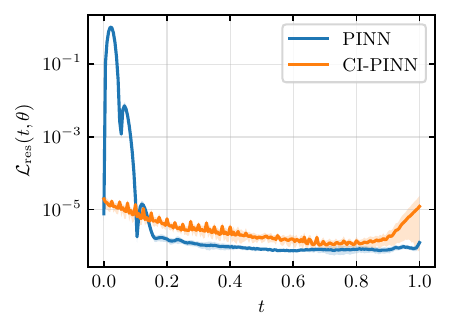}
    \caption{Residual loss}
    \label{fig:ac_loss}
  \end{subfigure}
  \hspace{10pt}
  \begin{subfigure}{0.45\textwidth}
    \includegraphics[width=\linewidth]{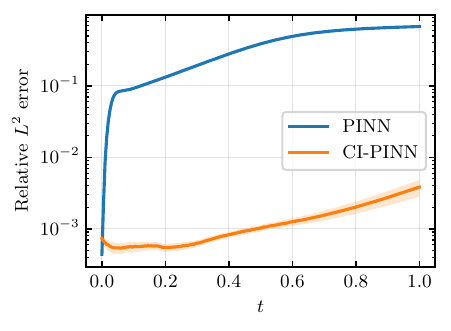}
    \caption{Relative $L^2$ error}
    \label{fig:ac_rl2e}
  \end{subfigure}
  \caption{Allen--Cahn equation. Time-dependent residual loss and relative $L^2$ error for PINN and CI-PINN. Both methods are trained on the same fixed grid with $(N_t,N_x)=(20,256)$. The CI-PINN result is included as a preview, and the method is introduced in Section~\ref{sec:methodology}.}
  \label{fig:ac_metrics}
\end{figure}

As shown in Figure~\ref{fig:ac_metrics}\subref{fig:ac_loss}, the standard PINN exhibits a pronounced residual peak near the initial time, followed by a rapid increase in the relative $L^2$ error in Figure~\ref{fig:ac_metrics}\subref{fig:ac_rl2e}. Although its residual decreases at later times, the solution error remains large. This behavior is consistent with the time-weighted error accumulation described by~\eqref{eq:error_propagation}. In comparison, CI-PINN produces a smaller early-time residual and maintains a substantially lower prediction error.

Taken together, the error estimate and the Allen--Cahn example show that the temporal distribution of the residual is relevant to trajectory accuracy. This motivates the history-dependent neural representation introduced in the next section.

\section{Methodology}\label{sec:methodology}
In this section, we develop the causal integral neural representation underlying CI-PINN. We begin by introducing CinNet, which extends the standard MLP with history-dependent integral features. We then describe how these features are evaluated using a prescribed temporal weight and fixed-mesh quadrature. Finally, we embed CinNet into the standard physics-informed training framework.

\subsection{Causal integral neural network}\label{sec:cinnet}
The integral representation in~\eqref{eq:evolution_integral_form} suggests that the preceding evolution should enter the network representation explicitly. Based on this observation, we introduce the causal integral neural network (CinNet), which augments each layer with a history-dependent integral pathway while retaining a direct pathway at the current coordinates.

Set $\mathbf{y}^{(0)}(t,\mathbf{x})=(t,\mathbf{x})$. For $\ell=1,\ldots,L$, we define the local and historical features by
\begin{align*}
  \mathbf{F}^{(\ell)}(t,\mathbf{x}) & = \mathbf{y}^{(\ell-1)}(t,\mathbf{x})\mathbf{W}_F^{(\ell)} + \mathbf{b}_F^{(\ell)}, \\
  \mathbf{V}^{(\ell)}(s,\mathbf{x}) & = \mathbf{y}^{(\ell-1)}(s,\mathbf{x})\mathbf{W}_V^{(\ell)} + \mathbf{b}_V^{(\ell)}.
\end{align*}
Here, $\mathbf{W}_F^{(\ell)}$ and $\mathbf{W}_V^{(\ell)}$ are trainable weight matrices, while $\mathbf{b}_F^{(\ell)}$ and $\mathbf{b}_V^{(\ell)}$ are the corresponding bias vectors. The local feature $\mathbf{F}^{(\ell)}$ transforms the representation at the current coordinates $(t,\mathbf{x})$, whereas the historical feature $\mathbf{V}^{(\ell)}$ maps representations at earlier times into the feature space used for temporal aggregation. The historical information available at time $t$ is then represented by the causal integral
\begin{equation}\label{eq:integral_feature}
  \mathbf{I}^{(\ell)}(t,\mathbf{x}) = \int_{t_0}^{t} \mathcal{A}^{(\ell)}(t,s)\mathbf{V}^{(\ell)}(s,\mathbf{x})\,\mathrm{d}s,
\end{equation}
where $\mathcal{A}^{(\ell)}(t,s)$ is a scalar temporal weight that controls the contribution of the feature evaluated at time $s$ to the representation at time $t$. Since the integration is restricted to $[t_0,t]$, the resulting integral feature depends only on information available up to the current time.

The local and integral features carry complementary information. Rather than combining them with fixed coefficients, CinNet introduces the learnable gating vector
\begin{equation*}
  \mathbf{z}^{(\ell)}=\operatorname{sigmoid}\left(\boldsymbol{\eta}^{(\ell)}\right),
\end{equation*}
where $\boldsymbol{\eta}^{(\ell)}$ is a trainable parameter vector. Each component of $\mathbf{z}^{(\ell)}$ therefore lies in $(0,1)$ and controls the balance between the corresponding channels of the local and integral features.

The layer output is then defined by
\begin{equation}\label{eq:cinnet_layer}
  \mathbf{y}^{(\ell)}(t,\mathbf{x}) = \sigma_\ell\left(\mathbf{z}^{(\ell)}\odot\mathbf{F}^{(\ell)}(t,\mathbf{x}) + \left(1-\mathbf{z}^{(\ell)}\right)\odot\mathbf{I}^{(\ell)}(t,\mathbf{x})\right).
\end{equation}
Here, $\odot$ denotes componentwise multiplication, and $\sigma_\ell$ acts componentwise, with $\sigma_L$ the identity. For $\ell>1$, $\mathbf{V}^{(\ell)}(s,\mathbf{x})$ is built from the already aggregated representation $\mathbf{y}^{(\ell-1)}(s,\mathbf{x})$, propagating history dependence through successive layers.

\begin{figure}[!htbp]
  \centering
  \includegraphics[width=0.9\linewidth]{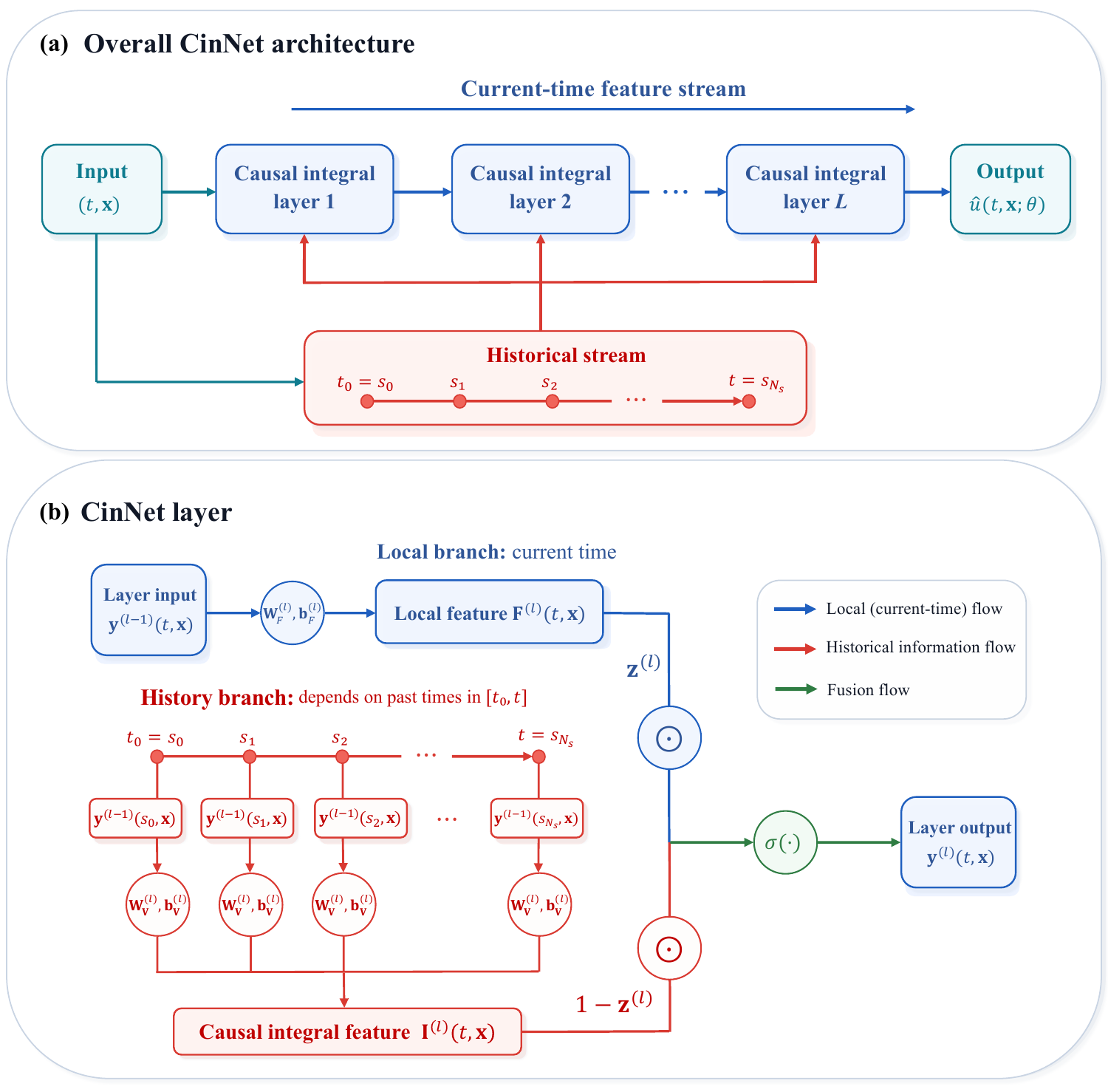}
  \caption{Schematic of CinNet. (a) Overall propagation of the current-time and historical feature streams through successive causal integral layers. (b) Construction of a single CinNet layer, in which the local feature and the history-dependent integral feature are combined through a learnable gate.}
  \label{fig:cinnet}
\end{figure}

The CinNet output is $\hat{u}(t,\mathbf{x};\theta)=\mathbf{y}^{(L)}(t,\mathbf{x})$, where $\theta$ collects all trainable parameters. Its causal integral~\eqref{eq:integral_feature} accumulates historical features in analogy with~\eqref{eq:evolution_integral_form}. The whole architecture is illustrated in Figure~\ref{fig:cinnet}.

\subsection{Practical implementation}\label{sec:implementation}
In this subsection, we describe the practical implementation of the causal integral term in CinNet, focusing on the choice of the temporal weight $\mathcal{A}^{(\ell)}(t,s)$ and the numerical approximation of the integral feature $\mathbf{I}^{(\ell)}(t,\mathbf{x})$.

\paragraph{Choice of the temporal weight.}
The general form~\eqref{eq:integral_feature} allows flexible choices of the temporal weight $\mathcal{A}^{(\ell)}(t,s)$. Although one may use a learnable weighting function depending on $(t,s,\mathbf{x})$, in this work we adopt a simple fixed temporal weight based on the temporal distance $t-s$. Specifically, for $t_0\le s\le t\le T$, we use the monomial form
\begin{equation}\label{eq:monomial_attention}
  \mathcal{A}^{(\ell)}(t,s) = \left(\frac{t-s}{T-t_0}\right)^q, \qquad q=2, \quad \forall \ell.
\end{equation}
The temporal weight~\eqref{eq:monomial_attention} is independent of $\mathbf{x}$ and introduces no additional trainable parameters. It assigns greater weight to features farther in the past, with the aim of preserving their contribution during temporal propagation. The exponent $q$ controls this preference: $q=0$ gives a constant weight, whereas any $q>0$ yields a weight that increases with the temporal distance $t-s$. We use $q=2$ as a fixed architectural choice throughout the numerical experiments.

\paragraph{Fixed-mesh quadrature.}
The integration interval in $\mathbf{I}^{(\ell)}(t,\mathbf{x})$ varies with the query time $t$. To evaluate this integral using a common set of temporal nodes, we introduce a uniform mesh on $[t_0,T]$,
\begin{equation*}
  t_0=s_0<s_1<\cdots<s_{N_s}=T,
\end{equation*}
where $s_{k+1}-s_k=\Delta s$ for $k=0,\ldots,N_s-1$.
For a query time $t\in(s_m,s_{m+1}]$, the first $m$ mesh intervals are fully contained in $[t_0,t]$, while $[s_m,t]$ provides the final contribution. We therefore define the effective length of the $k$th mesh interval by
\begin{equation*}
  r_k(t)=\max(0,t-s_k)-\max(0,t-s_{k+1})=\operatorname{clip}(t-s_k,0,\Delta s),
\end{equation*}
where $\operatorname{clip}(a,0,\Delta s)=\min\{\max\{a,0\},\Delta s\}$. Consequently, $r_k(t)=\Delta s$ for $k<m$, $r_m(t)=t-s_m$, and $r_k(t)=0$ for $k>m$.
Approximating the integrand on each active mesh interval by its value at the left endpoint gives
\begin{equation}\label{eq:left_rectangle}
  \mathbf{I}^{(\ell)}(t,\mathbf{x}) \approx \sum_{k=0}^{N_s-1} r_k(t)\mathcal{A}^{(\ell)}(t,s_k)\mathbf{V}^{(\ell)}(s_k,\mathbf{x}).
\end{equation}
Here, $r_k(t)$ specifies the portion of the $k$th mesh interval assigned to its left endpoint $s_k$, while $\mathcal{A}^{(\ell)}(t,s_k)$ weights the value $\mathbf{V}^{(\ell)}(s_k,\mathbf{x})$ in the quadrature sum.

The quadrature rule is incorporated into the network architecture for constructing historical features, rather than for obtaining a high-order approximation of a prescribed physical integral. This is also one reason for adopting the left-rectangle rule in our implementation. The quadrature node locations determine the past time instances at which historical features are sampled, whereas the node count $N_s$ controls the temporal resolution of the history representation.

The fixed mesh is particularly useful across layers. Computing $\mathbf{V}^{(\ell)}(s_k,\mathbf{x})$ requires $\mathbf{y}^{(\ell-1)}(s_k,\mathbf{x})$, whose integral term involves only the nodes $s_j\leq s_k$. All such integrals can therefore be evaluated on the same temporal mesh by restricting the quadrature sum to nodes no later than $s_k$. If a separate set of quadrature nodes were generated for every upper limit, each layer could introduce additional evaluation points required by the preceding layer, causing the number of nodes to grow with the network depth. The fixed mesh avoids this growth and provides a common quadrature structure across layers.

To evaluate~\eqref{eq:left_rectangle} at the $\ell$th layer, CinNet evaluates the preceding-layer output $\mathbf{y}^{(\ell-1)}$ at both the query point $(t,\mathbf{x})$ and the fixed quadrature points $\{(s_k,\mathbf{x})\}_{k=0}^{N_s-1}$. The value $\mathbf{y}^{(\ell-1)}(t,\mathbf{x})$ is used to compute the local term $\mathbf{F}^{(\ell)}(t,\mathbf{x})$, while the values $\mathbf{y}^{(\ell-1)}(s_k,\mathbf{x})$ are used to compute $\mathbf{V}^{(\ell)}(s_k,\mathbf{x})$ in the quadrature sum. These evaluations are performed at every layer. The fixed temporal nodes are used only for numerical quadrature and do not define an autoregressive time-stepping procedure. The surrogate $\hat{u}(t,\mathbf{x};\theta)$ can therefore be evaluated directly at arbitrary query times.

\subsection{Causal integral physics-informed neural network}\label{sec:cipinn_formulation}
In this subsection, we formulate CI-PINN for the evolution equation~\eqref{eq:pde}. We represent the approximate solution using CinNet,
\begin{equation*}
  \hat{u}(t,\mathbf{x};\theta) = \operatorname{CinNet}(t,\mathbf{x};\theta),
\end{equation*}
and determine the network parameters by minimizing the physics-informed objective in~\eqref{eq:loss}. Initial and boundary conditions are imposed through the corresponding loss terms. In the experiments below, periodic boundary conditions are imposed exactly by the input embedding in Appendix~\ref{app:exact_periodic_bc}, so the boundary loss is identically zero.

Rather than processing each space--time coordinate independently, CI-PINN constructs the solution representation at time $t$ by combining the current input with hidden representations from preceding times. For each query point $(t,\mathbf{x})$, a CinNet layer evaluates the representation from the preceding layer at the fixed temporal nodes $\{(s_k,\mathbf{x})\}_{k=0}^{N_s-1}$. These values define the causal integral term in~\eqref{eq:left_rectangle}, which is combined with the local term through~\eqref{eq:cinnet_layer}. Repeating this update across the $L$ layers gives $\hat{u}(t,\mathbf{x};\theta)$.

These intermediate evaluations involve only hidden representations and do not introduce additional residual collocation constraints. The residual collocation points and the initial- and boundary-condition points are used only in their corresponding loss terms. The complete procedure is summarized in Algorithm~\ref{alg:cipinn_training}.

\begin{algorithm}[htbp]
  \caption{CI-PINN}
  \label{alg:cipinn_training}
  \begin{algorithmic}[1]
    \Require Residual collocation set $\mathcal{D}_{\mathrm{res}}$, initial- and boundary-condition point sets, fixed temporal mesh $\{s_k\}_{k=0}^{N_s}$, and number of iterations $N_{\mathrm{iter}}$
    \Ensure Approximate solution $\hat{u}(t,\mathbf{x};\theta^\star)$
    \State Initialize the CinNet parameters $\theta$
    \For{$n=0,\ldots,N_{\mathrm{iter}}-1$}
    \State Evaluate $\hat{u}$ at the residual collocation points and the initial- and boundary-condition points using~\eqref{eq:cinnet_layer} and~\eqref{eq:left_rectangle}
    \State Compute the required derivatives by automatic differentiation and form the loss using~\eqref{eq:loss} and~\eqref{eq:loss_terms}
    \State Update $\theta$ using a gradient-based optimizer
    \EndFor
    \State Set $\theta^\star=\theta$
    \State \Return $\hat{u}(t,\mathbf{x};\theta^\star)$
  \end{algorithmic}
\end{algorithm}

\section{Numerical results}
\label{sec:results}

In this section, we present a series of numerical experiments to comprehensively evaluate the performance of the proposed CI-PINN for nonlinear evolution equations. The standard PINN is included as a canonical reference baseline, allowing us to quantify the improvement achieved by explicitly incorporating temporal causality into the network architecture. To provide a more competitive evaluation, we further compare CI-PINN with Causal PINN \cite{wang2024Respecting}, a representative causality-aware method that enforces temporal causality at the training level through adaptive temporal weighting of the residual loss. A detailed description of Causal PINN is provided in Appendix~\ref{app:causal_pinn}. This comparison enables us to assess the proposed architecture-level causal representation against an established training-level strategy for time-dependent problems. Through these experiments, we investigate predictive accuracy, training robustness, the effect of training-set size, sensitivity to selected hyperparameters, and the contributions of the main components of CinNet.

We first summarize the common computational settings used throughout the numerical experiments. Following~\cite{wang2024Respecting}, for the one-dimensional problems, the residual collocation points and initial-condition points are selected from fixed uniform grids and remain unchanged throughout training. For the two-dimensional problem, the temporal points are uniformly spaced, while the spatial points are generated by Latin hypercube sampling, as specified in Section~\ref{sec:ac2d}. Unless otherwise stated, neither adaptive resampling nor a time-marching strategy is employed. The reference solutions for the PDE test problems are generated using the Chebfun package, with a Fourier spectral discretization using 512 modes in space and the fourth-order exponential time-differencing Runge--Kutta scheme (ETDRK4) with a time step of $10^{-5}$. All experiments are implemented in JAX and conducted on a workstation equipped with an NVIDIA Tesla V100 GPU.

Regarding the neural network architectures, all methods employ networks with four hidden layers, each containing 128 neurons, and use the $\tanh$ activation function. For CI-PINN, the default number of quadrature nodes is set to $N_s=4N_t$, except in the experiments specifically examining the effect of $N_s\in\{4N_t,5N_t,6N_t,7N_t,8N_t\}$. All periodic boundary conditions are imposed exactly using the periodic feature embedding described in Appendix~\ref{app:exact_periodic_bc}. For two-dimensional problems, we set $M_x=M_y=2$. For one-dimensional problems, we set $M=5$ for fourth-order equations and $M=10$ for all other equations.

Unless otherwise specified, all models are trained for 300,000 epochs using the Adam optimizer. The initial learning rate is set to $10^{-3}$ and decayed by a factor of $0.9$ every 5,000 epochs. The loss weights in~\eqref{eq:loss} are set to $\omega_{\mathrm{res}}=1$ and $\omega_{\mathrm{ic}}=100$, which we found to provide stable and comparable training across all methods. To assess run-to-run variability, each experimental configuration is independently repeated five times using different random seeds. Unless otherwise specified, all quantitative results are reported as the mean and standard deviation over these five runs. For Causal PINN, unless the effect of the causality parameter $\epsilon$ in~\eqref{eq:temporal_weights} is specifically investigated, we report the best result obtained over $\epsilon\in\{10^{-4},10^{-3},10^{-2},10^{-1},10^{0},10^{1},10^{2},10^{3},10^{4}\}$.

For all experiments, prediction accuracy is quantified using the relative $L^2$ error (RL2E) evaluated on the test set:
\begin{equation} \label{errordef}
  \mathrm{RL2E} = \frac{ \left( \sum_{i=1}^{N_{\mathrm{test}}} \left\| u(t_i,\mathbf{x}_i)-\hat{u}(t_i,\mathbf{x}_i;\theta) \right\|_2^2 \right)^{1/2} }{ \left( \sum_{i=1}^{N_{\mathrm{test}}} \left\| u(t_i,\mathbf{x}_i) \right\|_2^2 \right)^{1/2} },
\end{equation}
where $u(t,\mathbf{x})$ is the reference solution and $\hat{u}(t,\mathbf{x};\theta)$ denotes the solution predicted by the corresponding method.

\subsection{Allen--Cahn equation}
\label{sec:ac}

We first consider the one-dimensional Allen--Cahn problem used as an illustrative example in~\cite{wang2024Respecting}:
\begin{equation}\label{eq:ac}
  \begin{alignedat}{2}
     & u_t-0.0001 u_{xx}+5(u^3-u) = 0, \qquad &  & (t,x)\in [0,1]\times[-1,1], \\
     & u(0,x) = x^2\cos(\pi x), \qquad        &  & x\in[-1,1],                 \\
     & u(t,-1) = u(t,1), \qquad               &  & t\in[0,1].
  \end{alignedat}
\end{equation}
This problem is challenging for the standard continuous-time PINN formulation of Raissi et al.~\cite{raissi2019Physicsinformed}. To improve the resulting PINN predictions, Wight and Zhao~\cite{wight2021Solving} and McClenny and Braga-Neto~\cite{mcclenny2023Selfadaptive} proposed adaptive resampling and self-adaptive weighting strategies, respectively.

\begin{figure}[!htbp]
  \centering
  \begin{subfigure}{0.45\textwidth}
    \includegraphics[width=\linewidth]{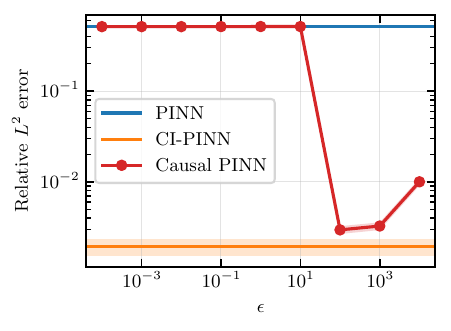}
    \caption{Causal PINN: causality parameter $\epsilon$}
    \label{fig:ac_eps}
  \end{subfigure}
  \hspace{10pt}
  \begin{subfigure}{0.45\textwidth}
    \includegraphics[width=\linewidth]{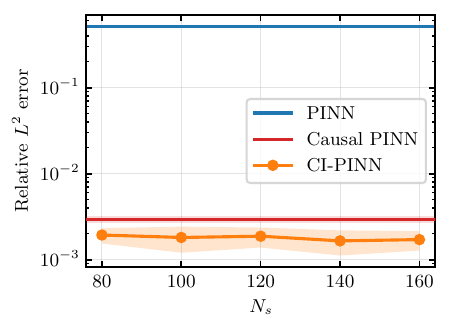}
    \caption{CI-PINN: number of quadrature nodes $N_s$}
    \label{fig:ac_ns}
  \end{subfigure}
  \caption{Allen--Cahn equation. Relative $L^2$ errors~\eqref{errordef} versus the causality parameter $\epsilon$ in~\eqref{eq:temporal_weights} for Causal PINN and the number of quadrature nodes $N_s$ in~\eqref{eq:left_rectangle} for CI-PINN.}
  \label{fig:ac_hyper_params}
\end{figure}

We first examine the sensitivity of Causal PINN and CI-PINN to their respective method-specific hyperparameters, as shown in Figure~\ref{fig:ac_hyper_params}. For Causal PINN, the additional hyperparameter of interest is the causality parameter $\epsilon$ in~\eqref{eq:temporal_weights}, which controls the strength of temporal weighting in the residual loss. For CI-PINN, the additional hyperparameter of interest is the number of quadrature nodes $N_s$ in~\eqref{eq:left_rectangle}, which determines the temporal resolution of the historical information available to the network. In both cases, we vary one hyperparameter while keeping all other hyperparameters fixed, and both methods are trained using the same fixed training set with $(N_t,N_x)=(20,256)$.

Figure~\ref{fig:ac_hyper_params} reveals a marked difference in the
hyperparameter sensitivity of the two methods. Causal PINN is highly
sensitive to the causality parameter $\epsilon$ and therefore requires
careful tuning. Without an appropriate choice of $\epsilon$, its performance
can fall back to a level comparable to that of the standard PINN.

In contrast, the performance of CI-PINN becomes stable once the number of
quadrature nodes reaches a moderate level. In the present example, the error
shows only minor variation for $N_s\geq 4N_t$, indicating that $N_s$ does
not require delicate tuning once the temporal history is resolved with
sufficient quadrature resolution. Moreover, CI-PINN consistently outperforms
both the standard PINN and Causal PINN with its best-performing tested value
$\epsilon=10^{2}$. Thus, although both methods introduce an additional
method-specific hyperparameter, CI-PINN incurs a substantially smaller
hyperparameter-tuning burden. Accordingly, we set $N_s=4N_t$ for CI-PINN and $\epsilon=10^{2}$ for Causal PINN in the remainder of this example.

We next investigate the influence of training grid resolution by varying $N_t$ and $N_x$. The relative $L^2$ errors are reported in Table~\ref{tab:ac_rl2e_NtNx}.

\begin{table}[!htbp]
  \centering
  \caption{Allen--Cahn equation. Relative $L^2$ errors of different methods for different training grid sizes $(N_t,N_x)$. The reported values are the mean and standard deviation over five independent runs. For Causal PINN, the selected parameter $\epsilon=10^2$ is used.}
  \begin{tabular}{ccccc}
    \toprule
    $N_t$ & $N_x$                                     & PINN & Causal PINN & CI-PINN \\
    \midrule
    \multirow{3}{*}{10}
          & 64
          & $(6.98\pm0.17)\times10^{-1}$
          & $(2.45\pm0.32)\times10^{-1}$
          & $\boldsymbol{(1.35\pm0.54)\times10^{-2}}$                                \\
          & 128
          & $(5.16\pm0.04)\times10^{-1}$
          & $(1.98\pm1.29)\times10^{-2}$
          & $\boldsymbol{(2.99\pm0.31)\times10^{-3}}$                                \\
          & 256
          & $(5.16\pm0.03)\times10^{-1}$
          & $(1.13\pm0.33)\times10^{-2}$
          & $\boldsymbol{(3.96\pm0.53)\times10^{-3}}$                                \\
    \midrule
    \multirow{3}{*}{20}
          & 64
          & $(6.96\pm0.19)\times10^{-1}$
          & $(9.40\pm4.10)\times10^{-2}$
          & $\boldsymbol{(1.35\pm0.51)\times10^{-2}}$                                \\
          & 128
          & $(5.13\pm0.01)\times10^{-1}$
          & $(4.15\pm1.34)\times10^{-3}$
          & $\boldsymbol{(1.80\pm0.48)\times10^{-3}}$                                \\
          & 256
          & $(5.15\pm0.01)\times10^{-1}$
          & $(2.96\pm0.29)\times10^{-3}$
          & $\boldsymbol{(1.94\pm0.39)\times10^{-3}}$                                \\
    \bottomrule
  \end{tabular}
  \label{tab:ac_rl2e_NtNx}
\end{table}

As shown in Table~\ref{tab:ac_rl2e_NtNx}, the standard PINN remains inaccurate across all training grid sizes considered. While Causal PINN reduces the error to varying degrees, CI-PINN consistently achieves the smallest relative $L^2$ error for every tested $(N_t,N_x)$ pair. This consistent performance advantage demonstrates the efficacy of incorporating temporal causality at the representation level, rather than solely at the training level as implemented in the baseline methods.

The advantage of CI-PINN is particularly pronounced in the sparse-collocation regime. For example, when $(N_t,N_x)=(10,64)$, the baseline methods still produce dynamically inaccurate solutions with large errors, whereas CI-PINN gives a reasonably accurate prediction. This behavior suggests an implicit regularization effect of the causal integral representation under sparse collocation.

To further assess collocation efficiency, we perform an additional dense-grid experiment with the standard PINN. Using $(N_t,N_x)=(100,256)$, the standard PINN achieves a relative $L^2$ error of $4.17\times10^{-3}$. In comparison, CI-PINN achieves a slightly smaller mean error of $2.99\times10^{-3}$ using only $(N_t,N_x)=(10,128)$, which corresponds to a 20-fold reduction in the number of residual collocation points. This comparison demonstrates that CI-PINN makes more effective use of sparse space--time collocation points.

Finally, we investigate the sensitivity of the three methods to the initial-condition loss weight $\omega_{\mathrm{ic}}$. This parameter is important in PINN training because it controls the balance between fitting the initial condition and minimizing the PDE residual. In this experiment, we fix the training grid at $(N_t,N_x)=(20,256)$ and vary $\omega_{\mathrm{ic}}\in\{10^0, 10^1, 10^2, 10^3, 10^4\}$. The results are reported in Table~\ref{tab:ac_rl2e_wi}.

\begin{table}[!htbp]
  \centering
  \caption{Allen--Cahn equation. Relative $L^2$ errors of different methods for different initial-condition loss weights $\omega_{\mathrm{ic}}$. The reported values are the mean and standard deviation over five independent runs. For Causal PINN, the selected parameter $\epsilon=10^2$ is used.}
  \begin{tabular}{cccc}
    \toprule
    $\omega_{\mathrm{ic}}$ & PINN                                      & Causal PINN & CI-PINN \\
    \midrule
    $10^0$
                           & $(6.93\pm3.82)\times10^{-1}$
                           & $(3.94\pm1.85)\times10^{-3}$
                           & $\boldsymbol{(1.51\pm0.36)\times10^{-3}}$                         \\
    $10^1$
                           & $(4.79\pm0.45)\times10^{-1}$
                           & $(1.35\pm0.44)\times10^{-3}$
                           & $\boldsymbol{(7.05\pm1.56)\times10^{-4}}$                         \\
    $10^2$
                           & $(5.15\pm0.01)\times10^{-1}$
                           & $(2.96\pm0.29)\times10^{-3}$
                           & $\boldsymbol{(1.94\pm0.39)\times10^{-3}}$                         \\
    $10^3$
                           & $(5.16\pm0.01)\times10^{-1}$
                           & $(6.35\pm1.28)\times10^{-3}$
                           & $\boldsymbol{(3.07\pm0.75)\times10^{-3}}$                         \\
    $10^4$
                           & $(5.49\pm0.03)\times10^{-1}$
                           & $(4.33\pm0.26)\times10^{-1}$
                           & $\boldsymbol{(6.54\pm1.95)\times10^{-3}}$                         \\
    \bottomrule
  \end{tabular}
  \label{tab:ac_rl2e_wi}
\end{table}

As shown in Table~\ref{tab:ac_rl2e_wi}, CI-PINN yields the smallest relative $L^2$ error for every tested value of $\omega_{\mathrm{ic}}$ and maintains accuracy consistently across the entire range. Causal PINN performs well for moderate values of $\omega_{\mathrm{ic}}$, but its error increases sharply when $\omega_{\mathrm{ic}}=10^4$. The standard PINN remains inaccurate for all tested weights. These results indicate that CI-PINN is more robust to the choice of $\omega_{\mathrm{ic}}$, which is a desirable property in practical applications because the optimal loss weighting may not be known a priori.

\subsection{Korteweg--de Vries equation}
\label{sec:kdv}

We next consider the Korteweg--de Vries (KdV) equation as a nonlinear dispersive evolution problem. Unlike the dissipative dynamics of the Allen--Cahn equation, the KdV equation describes nonlinear wave propagation governed by a balance between nonlinearity and dispersion. The problem is given by
\begin{equation}\label{eq:kdv}
  \begin{alignedat}{2}
     & u_t+\eta uu_x+\mu^2 u_{xxx} = 0, \qquad &  & (t,x)\in [0,1]\times[-1,1], \\
     & u(0,x) = \cos(\pi x), \qquad            &  & x\in[-1,1],                 \\
     & u(t,-1) = u(t,1), \qquad                &  & t\in[0,1].
  \end{alignedat}
\end{equation}
Here, $\eta$ controls the strength of the nonlinear term, whereas $\mu$
determines the strength of dispersion. Under the KdV dynamics \eqref{eq:kdv}, the prescribed
initial profile evolves into a train of solitary-wave structures. Following
the parameter setting in \cite{zabusky1965Interaction}, we set $\eta=1$ and
$\mu=0.022$.

\begin{figure}[!htbp]
  \centering
  \includegraphics[width=0.9\linewidth]{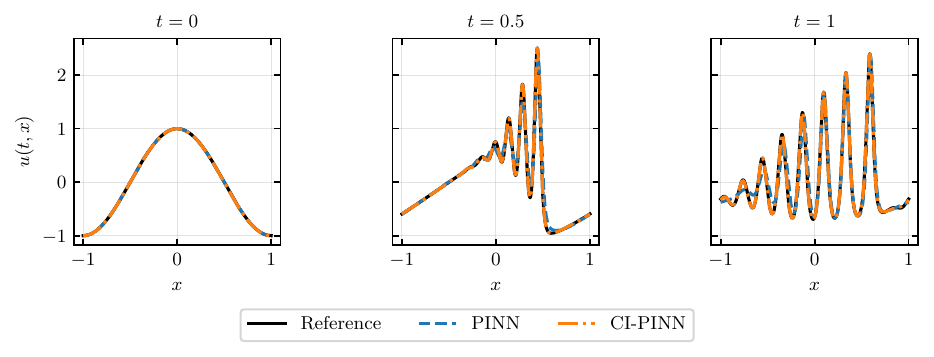}
  \caption{KdV equation. Comparison of the reference solution with the PINN and CI-PINN solution profiles at $t=0$, $0.5$, and $1$. For each method, the prediction from the best-performing run among five independent runs is shown. The corresponding relative $L^2$ errors are $1.87\times10^{-1}$ for PINN and $6.94\times10^{-3}$ for CI-PINN.}
  \label{fig:kdv_snapshots}
\end{figure}

We first compare the reference solution of the KdV equation with the predictions obtained by PINN and CI-PINN, as shown in Figure~\ref{fig:kdv_snapshots}. Both methods are trained on the same fixed space--time grid with $(N_t,N_x)=(20,256)$, and CI-PINN uses $N_s=4N_t=80$ quadrature nodes.

As shown in Figure~\ref{fig:kdv_snapshots}, the solution develops increasingly complex oscillatory wave structures as time evolves. Although PINN captures the overall evolution pattern, noticeable discrepancies appear in the amplitudes and phases of the oscillations, particularly at later times. These local errors accumulate across the rapidly varying wave structures, resulting in a substantial deviation from the reference solution. In contrast, CI-PINN closely follows the reference profiles and accurately captures both the locations and amplitudes of the oscillatory structures throughout the evolution. Over five independent runs, the mean relative $L^2$ error is reduced from $(4.74\pm2.76)\times10^{-1}$ for PINN to $(7.81\pm0.68)\times10^{-3}$ for CI-PINN. These results demonstrate that the causal integral representation substantially improves the approximation of the nonlinear dispersive dynamics under the same limited space--time collocation setting.

\begin{figure}[!htbp]
  \centering
  \begin{subfigure}{0.45\textwidth}
    \includegraphics[width=\linewidth]{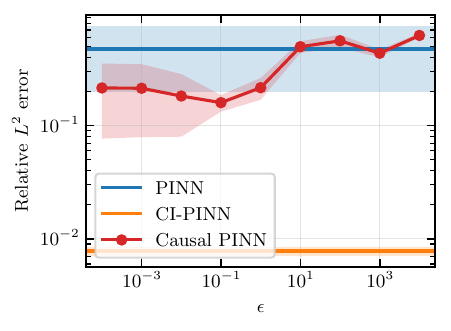}
    \caption{Causal PINN: causality parameter $\epsilon$}
    \label{fig:kdv_eps}
  \end{subfigure}
  \hspace{10pt}
  \begin{subfigure}{0.45\textwidth}
    \includegraphics[width=\linewidth]{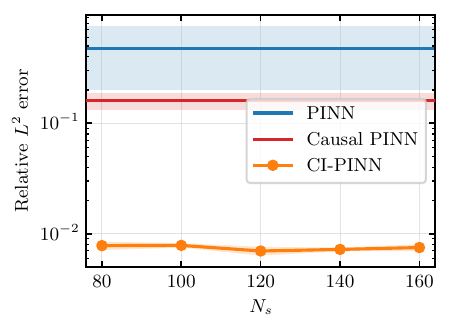}
    \caption{CI-PINN: number of quadrature nodes $N_s$}
    \label{fig:kdv_ns}
  \end{subfigure}
  \caption{KdV equation. Relative $L^2$ errors~\eqref{errordef} versus the causality parameter $\epsilon$ in~\eqref{eq:temporal_weights} for Causal PINN and the number of quadrature nodes $N_s$ in~\eqref{eq:left_rectangle} for CI-PINN.}
  \label{fig:kdv_hyper_params}
\end{figure}

Figure~\ref{fig:kdv_hyper_params} compares the sensitivity of Causal PINN
and CI-PINN to their respective method-specific hyperparameters for the
KdV equation. All methods are trained on the same fixed training grid with
$(N_t,N_x)=(20,256)$, and the reported errors are the mean and standard
deviation over five independent runs.

As shown in Figure~\ref{fig:kdv_hyper_params}\subref{fig:kdv_eps}, Causal PINN yields relatively large errors for all the tested values of the causality parameter $\epsilon$, and its performance remains broadly comparable to that of the standard PINN. Even at its best-performing tested value, $\epsilon=10^{-1}$, Causal PINN attains a relative $L^2$ error of $(1.60\pm0.27)\times10^{-1}$, whereas the standard PINN yields $(4.74\pm2.76)\times10^{-1}$. By comparison, CI-PINN with $N_s=80$ achieves a substantially smaller error of $(7.81\pm0.68)\times10^{-3}$. Thus, even after tuning $\epsilon$, Causal PINN provides only a limited improvement over the standard PINN, while CI-PINN improves the prediction accuracy by more than one order of magnitude relative to both baseline methods.

Consistent with the Allen--Cahn results in Figure~\ref{fig:ac_hyper_params}\subref{fig:ac_ns}, Figure~\ref{fig:kdv_hyper_params}\subref{fig:kdv_ns} shows that CI-PINN maintains stable accuracy over $N_s\in\{80,100,120,140,160\}$. This further confirms that no delicate tuning of $N_s$ is required once a moderate quadrature resolution is reached.

\subsection{Cahn--Hilliard equation}\label{sec:ch}

We next consider the Cahn--Hilliard equation, which introduces additional
challenges through the combination of nonlinear phase-separation dynamics
and a fourth-order spatial derivative. The resulting evolution involves
sharp interfacial structures and places greater demands on the approximation
of high-order nonlinear space--time dynamics. The governing equation is
\begin{equation}\label{eq:ch}
  \begin{alignedat}{2}
     & u_t+10^{-6}u_{xxxx}-10^{-2}(u^3-u)_{xx}=0, \qquad &  & (t,x)\in (0,1]\times[-1,1], \\
     & u(0,x) = -\cos(2\pi x), \qquad                    &  & x \in [-1,1],               \\
     & u(t,-1) = u(t,1), \qquad                          &  & t \in [0,1].
  \end{alignedat}
\end{equation}

\begin{figure}[!htbp]
  \centering
  \includegraphics[width=0.9\linewidth]{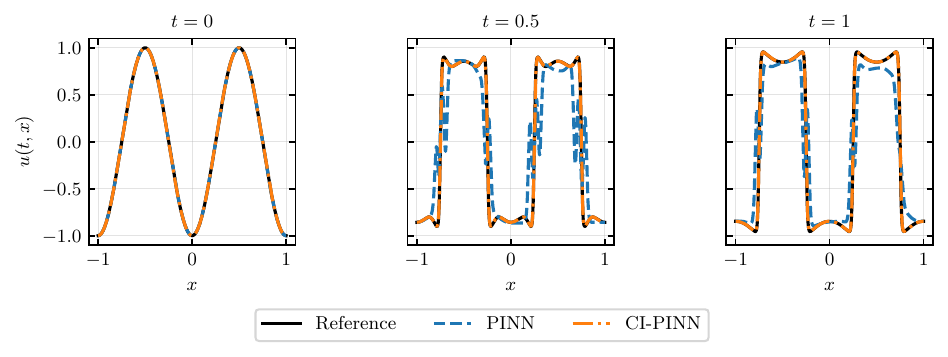}
  \caption{Cahn--Hilliard equation. Comparison of the reference solution with the PINN and CI-PINN solution profiles at $t=0$, $0.5$, and $1$. For each method, the prediction from the best-performing run among five independent runs is shown. The corresponding relative $L^2$ errors are $4.07\times10^{-1}$ for PINN and $5.68\times10^{-3}$ for CI-PINN.}
  \label{fig:ch_snapshots}
\end{figure}

For the Cahn--Hilliard problem~\eqref{eq:ch}, we first compare the reference solution with the predictions obtained by PINN and CI-PINN, as shown in Figure~\ref{fig:ch_snapshots}. Both methods are trained on the same fixed space--time grid with $(N_t,N_x)=(40,256)$, and CI-PINN uses $N_s=4N_t=160$ quadrature nodes.

As shown in Figure~\ref{fig:ch_snapshots}, PINN captures the overall solution pattern but exhibits pronounced discrepancies near the sharp interfaces, particularly at later times. In contrast, CI-PINN closely follows the reference profiles and substantially improves the resolution of the interfacial structures. Over five independent runs, the mean relative $L^2$ error is reduced from $(4.26\pm0.18)\times10^{-1}$ for PINN to $(6.06\pm0.28)\times10^{-3}$ for CI-PINN. These results indicate that the causal integral representation enables CI-PINN to capture the nonlinear phase-separation dynamics much more accurately under the same limited space--time collocation setting.

\begin{figure}[!htbp]
  \centering
  \begin{subfigure}{0.45\textwidth}
    \includegraphics[width=\linewidth]{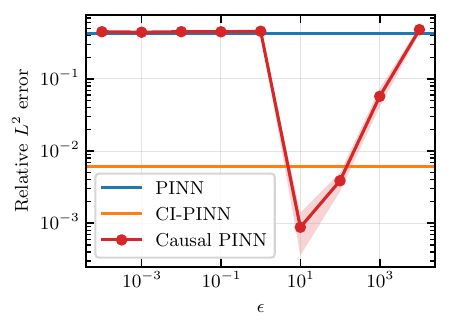}
    \caption{Causal PINN: causality parameter $\epsilon$}
    \label{fig:ch_eps}
  \end{subfigure}
  \hspace{10pt}
  \begin{subfigure}{0.45\textwidth}
    \includegraphics[width=\linewidth]{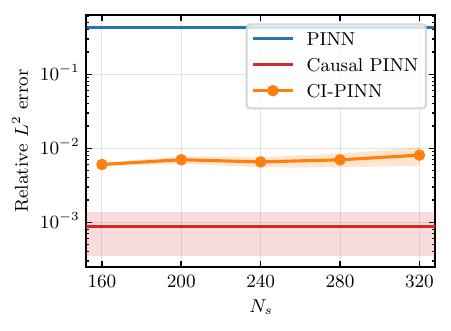}
    \caption{CI-PINN: number of quadrature nodes $N_s$}
    \label{fig:ch_ns}
  \end{subfigure}
  \caption{Cahn--Hilliard equation. Relative $L^2$ errors~\eqref{errordef} versus the causality parameter $\epsilon$ in~\eqref{eq:temporal_weights} for Causal PINN and the number of quadrature nodes $N_s$ in~\eqref{eq:left_rectangle} for CI-PINN.}
  \label{fig:ch_hyper_params}
\end{figure}

We next examine the dependence of Causal PINN and CI-PINN on their respective method-specific hyperparameters, as shown in Figure~\ref{fig:ch_hyper_params}. All methods are trained on the same fixed training grid with $(N_t,N_x)=(40,256)$.

As shown in Figure~\ref{fig:ch_hyper_params}\subref{fig:ch_eps}, Causal PINN is highly sensitive to the causality parameter $\epsilon$. For most of the tested values, its performance remains comparable to that of the standard PINN, which yields a relative $L^2$ error of $(4.26\pm0.18)\times10^{-1}$. A substantial improvement is obtained only within a relatively narrow range of $\epsilon$. At the best-performing tested value $\epsilon=10$, Causal PINN achieves a relative $L^2$ error of $(8.78\pm5.22)\times10^{-4}$. By comparison, CI-PINN with $N_s=160$ achieves an error of $(6.06\pm0.28)\times10^{-3}$. Thus, although Causal PINN with its best tested parameter attains the highest accuracy in this example, this performance relies strongly on a carefully selected causality parameter, whereas inappropriate choices of $\epsilon$ lead to a pronounced degradation in accuracy.

In contrast, Figure~\ref{fig:ch_hyper_params}\subref{fig:ch_ns} shows that CI-PINN maintains stable
accuracy across different values of $N_s$, and consistently achieves
errors on the order of $10^{-3}$, substantially below those of the standard
PINN. Consistent with the Allen--Cahn and KdV results in
Figures~\ref{fig:ac_hyper_params}\subref{fig:ac_ns} and \ref{fig:kdv_hyper_params}\subref{fig:kdv_ns}, only minor variations are
observed as $N_s$ changes. This further confirms that CI-PINN does not
require delicate tuning of the quadrature resolution once a moderate number
of quadrature nodes is used.

Finally, we conduct an ablation study to evaluate the contribution of the architectural components in CinNet for the Cahn--Hilliard equation. The results are summarized in Table~\ref{tab:ch_ablation}.

\begin{table}[!htbp]
  \centering
  \caption{Cahn--Hilliard equation. Ablation study of the architectural components in CinNet with $\omega_{\mathrm{ic}}=100$, $(N_t,N_x)=(40,256)$, and $N_s=160$. For the temporal weight, \xmark~denotes $\mathcal{A}^{(\ell)}(t,s)=1$. For the gating mechanism, \xmark~denotes that $\mathbf{z}^{(\ell)}$ and $1-\mathbf{z}^{(\ell)}$ are both fixed to 0.5 in~\eqref{eq:cinnet_layer}.}
  \begin{tabular}{ccc}
    \toprule
    Temporal weight & Gating mechanism & Relative $L^2$ error                       \\
    \midrule
    \cmark          & \cmark           & $\boldsymbol{(6.06\pm0.28)\times 10^{-3}}$ \\
    \cmark          & \xmark           & $(1.95\pm0.52)\times 10^{-1}$              \\
    \xmark          & \cmark           & $(3.35\pm0.39)\times 10^{-1}$              \\
    \xmark          & \xmark           & $(4.31\pm0.22)\times 10^{-1}$              \\
    \bottomrule
  \end{tabular}
  \label{tab:ch_ablation}
\end{table}

As shown in Table~\ref{tab:ch_ablation}, only the CI-PINN variant retaining both the temporal weight and gating mechanism achieves a relatively low relative $L^2$ error. Removing either component causes a substantial increase in the error, indicating that both components are essential for obtaining a reliable prediction of the Cahn--Hilliard dynamics.

\subsection{Two-dimensional Allen--Cahn equation}
\label{sec:ac2d}

As the final benchmark, we consider a two-dimensional Allen--Cahn equation
to further assess the capability of CI-PINN for nonlinear evolution problems
in higher-dimensional spatial domains. Compared with the one-dimensional
examples considered above, this problem involves substantially richer
spatial structures and requires the network to resolve the temporal evolution
of sharp interfaces over a two-dimensional domain. The governing equation is
\begin{equation}\label{eq:ac2d}
  \begin{alignedat}{2}
     & u_t-10^{-4} \Delta u+5(u^3-u) = 0, \qquad                   &  & (t,x,y)\in [0,1]\times[-1,1]^2, \\
     & u(0,x,y)=\cos(\pi x)\cos(\pi y)(1-\exp(-(x^2+y^2))), \qquad &  & (x,y)\in[-1,1]^2,               \\
     & u(t,-1,y) = u(t,1,y), \quad u(t,x,-1) = u(t,x,1), \qquad    &  & t\in[0,1], \quad x,y \in[-1,1].
  \end{alignedat}
\end{equation}

\begin{figure}[!htbp]
  \centering
  \begin{subfigure}{0.9\textwidth}
    \includegraphics[width=\textwidth]{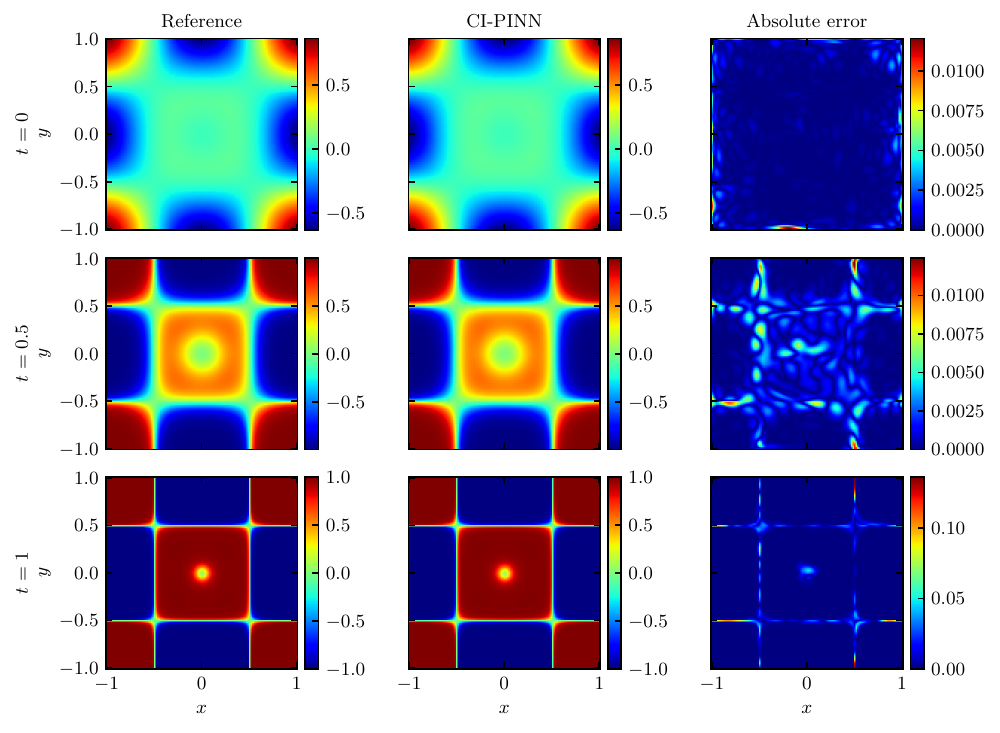}
    \caption{Spatial snapshots at $t=0,0.5,1.0$.}
    \label{fig:ac2d_heatmaps_t}
  \end{subfigure}

  \vspace{0.6em}

  \begin{subfigure}{0.9\textwidth}
    \includegraphics[width=\linewidth]{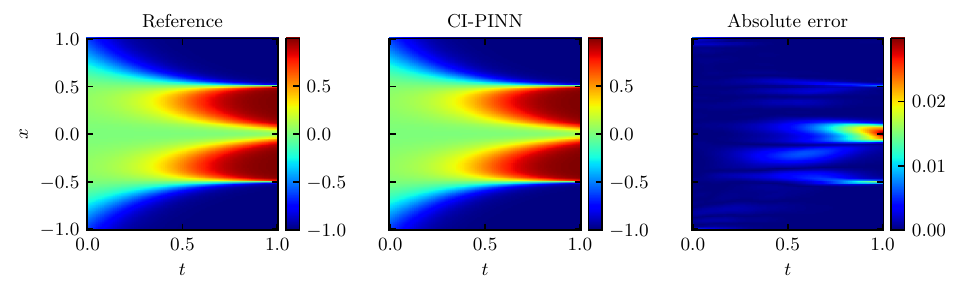}
    \caption{Space--time cross sections at $y=0$.}
    \label{fig:ac2d_snapshots_y}
  \end{subfigure}
  \caption{Two-dimensional Allen--Cahn equation. Comparison of the reference solution with the CI-PINN prediction and pointwise absolute error. CI-PINN is trained with $(N_t,N_{xy})=(20,1024)$ and uses $N_s=80$ quadrature nodes. Panel (a) shows the two-dimensional spatial snapshots at $t=0$, $0.5$, and $1.0$, where each row corresponds to one time level and the three columns show the reference solution, prediction, and pointwise absolute error, respectively. Panel (b) shows the corresponding space--time cross sections on the plane $y=0$, with $t$ on the horizontal axis and $x$ on the vertical axis. The relative $L^2$ error of CI-PINN is $4.29\times10^{-3}$ from the best-performing run among five independent runs.}
  \label{fig:ac2d_compare}
\end{figure}

We first examine the CI-PINN prediction for~\eqref{eq:ac2d}. The model is trained on a fixed space--time training set with $(N_t,N_{xy})=(20,1024)$, where $N_{xy}$ denotes the number of spatial collocation points. The temporal points are taken from a uniform grid, while a fixed set of $N_{xy}=1024$ spatial collocation points is generated by Latin hypercube sampling over $[-1,1]^2$ and shared across all temporal levels. CI-PINN uses $N_s=4N_t=80$ quadrature nodes. Figure~\ref{fig:ac2d_compare} compares the CI-PINN prediction with the reference solution and shows the pointwise absolute error. Figure~\ref{fig:ac2d_compare}\subref{fig:ac2d_heatmaps_t} shows two-dimensional spatial snapshots at $t=0$, $0.5$, and $1$, together with the corresponding pointwise absolute errors, whereas Figure~\ref{fig:ac2d_compare}\subref{fig:ac2d_snapshots_y} displays the space--time cross section on $y=0$.

As shown in Figure~\ref{fig:ac2d_compare}\subref{fig:ac2d_heatmaps_t}, the solution undergoes a pronounced nonlinear evolution from the smooth initial profile toward states characterized by increasingly sharp interfacial structures. CI-PINN accurately reproduces both the large-scale spatial pattern and the localized interfaces at all three displayed times. The pointwise errors remain small over most of the spatial domain and are primarily concentrated near the rapidly varying interfaces, where the solution is most difficult to approximate. The space--time cross sections in Figure~\ref{fig:ac2d_compare}\subref{fig:ac2d_snapshots_y} further show that CI-PINN closely tracks the temporal evolution of the reference solution throughout the interval $t\in[0,1]$. These results demonstrate that the causal integral representation remains effective when extended to a two-dimensional spatial domain and can accurately capture nonlinear interfacial dynamics under sparse space--time collocation.

\begin{figure}[!htbp]
  \centering
  \begin{subfigure}{0.45\textwidth}
    \includegraphics[width=\linewidth]{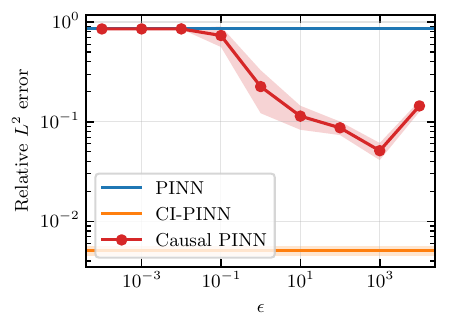}
    \caption{Causal PINN: causality parameter $\epsilon$}
    \label{fig:ac2d_eps}
  \end{subfigure}
  \hspace{10pt}
  \begin{subfigure}{0.45\textwidth}
    \includegraphics[width=\linewidth]{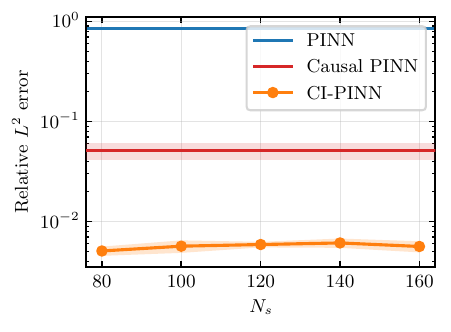}
    \caption{CI-PINN: number of quadrature nodes $N_s$}
    \label{fig:ac2d_ns}
  \end{subfigure}
  \caption{Two-dimensional Allen--Cahn equation. Relative $L^2$ errors~\eqref{errordef} versus the causality parameter $\epsilon$ in~\eqref{eq:temporal_weights} for Causal PINN and the quadrature node counts $N_s$ in~\eqref{eq:left_rectangle} for CI-PINN.}
  \label{fig:ac2d_hyper_params}
\end{figure}

We next investigate the dependence of Causal PINN and CI-PINN on their respective method-specific hyperparameters. All methods use the same fixed training set with $(N_t,N_{xy})=(20,1024)$. The results are shown in Figure~\ref{fig:ac2d_hyper_params}.

As shown in Figure~\ref{fig:ac2d_hyper_params}\subref{fig:ac2d_eps}, the performance of Causal
PINN depends strongly on the causality parameter $\epsilon$. For small values
of $\epsilon$, its accuracy remains close to that of the standard PINN,
whose relative $L^2$ error is $(8.53\pm0.03)\times10^{-1}$. A substantial
improvement is obtained only when $\epsilon$ is sufficiently large. Among
the tested values, the best performance is achieved at $\epsilon=10^{3}$,
for which Causal PINN attains a relative $L^2$ error of
$(5.12\pm1.00)\times10^{-2}$. By comparison, CI-PINN with $N_s=80$
achieves a markedly smaller error of $(5.09\pm0.54)\times10^{-3}$.
Thus, even after tuning the causality parameter, Causal PINN remains about
one order of magnitude less accurate than CI-PINN for this
two-dimensional problem.

Consistent with the one-dimensional results, Figure~\ref{fig:ac2d_hyper_params}\subref{fig:ac2d_ns} shows that CI-PINN remains insensitive to the choice of $N_s$ in the two-dimensional setting once $N_s\geq 4N_t$.

\section{Conclusions}\label{sec:conclusion}

We have presented CI-PINN, a causal integral physics-informed neural network for evolution equations, which introduces a Volterra-type causal integral term at the representation level to embed temporal causality directly into the solution representation. We have also provided a rationale for the proposed architecture and detailed its practical implementation. Numerical comparisons with standard PINN and Causal PINN demonstrate improvements in solution accuracy and training robustness, while the sensitivity studies show that the quadrature node count requires little tuning once a sufficient temporal resolution is reached. The sparse-collocation results further suggest an implicit regularization effect of the proposed representation. Systematic ablation studies verify the contributions of both the temporally weighted causal integral term and the gating mechanism. These results demonstrate the effectiveness of incorporating temporal causality directly into the neural representation for solving evolution equations.

Future work will investigate adaptive strategies for selecting quadrature nodes and temporal weights, with the aim of improving both efficiency and flexibility in history aggregation. Because the causal integral is taken only over time, extending the proposed architecture to higher-dimensional spatial problems beyond the one- and two-dimensional cases considered here is an important direction for future work. It would also be of interest to combine the proposed representation with alternative training strategies, adaptive sampling methods, and advanced optimization techniques, as well as to further characterize the theoretical effect of the causal integral representation on the approximation and optimization properties of physics-informed neural networks.

\section*{Code and data availability}
The code and data supporting the findings of this study are available from the authors upon reasonable request.

\section*{Appendix}
\appendix
\renewcommand{\thesection}{\Alph{section}}

\makeatletter
\@addtoreset{table}{section}
\@addtoreset{figure}{section}
\@addtoreset{theorem}{section}
\makeatother

\renewcommand{\thetable}{\thesection.\arabic{table}}
\renewcommand{\thefigure}{\thesection.\arabic{figure}}
\renewcommand{\thetheorem}{\thesection.\arabic{theorem}}

\section{Periodic feature embedding for exact boundary enforcement}\label{app:exact_periodic_bc}
Following~\cite{dong2021Method,wang2024Respecting}, periodic boundary conditions can be enforced exactly by using periodic feature embeddings for the spatial inputs. We first describe the one-dimensional case. For a spatial input $x\in[x_l,x_r]$, we replace $x$ with
\begin{equation}\label{eq:periodic_embedding_1d}
  v(x) = \left( 1,\cos(\omega_x x),\sin(\omega_x x), \ldots, \cos(M\omega_x x),\sin(M\omega_x x) \right),
\end{equation}
where $\omega_x=2\pi/(x_r-x_l)$ is the fundamental frequency and $M$ is the number of harmonics. The neural network then takes $(t,v(x))$ as input instead of $(t,x)$. Since the embedding~\eqref{eq:periodic_embedding_1d} is periodic with period $x_r-x_l$, the resulting network output is periodic in $x$ by construction.

For two-dimensional spatial domains, the periodic feature embedding is constructed by using tensor products of trigonometric features in the two coordinate directions. Let $(x,y)\in[x_l,x_r]\times[y_l,y_r]$. The two-dimensional periodic feature embedding is given by
\begin{equation}\label{eq:periodic_embedding_2d}
  v(x,y) = \left[
    \begin{array}{cccc}
      \cos(\omega_x x)\cos(\omega_y y), & \cos(\omega_x x)\cos(2\omega_y y), & \cdots & \cos(M_x\omega_x x)\cos(M_y\omega_y y) \\
      \cos(\omega_x x)\sin(\omega_y y), & \cos(\omega_x x)\sin(2\omega_y y), & \cdots & \cos(M_x\omega_x x)\sin(M_y\omega_y y) \\
      \sin(\omega_x x)\cos(\omega_y y), & \sin(\omega_x x)\cos(2\omega_y y), & \cdots & \sin(M_x\omega_x x)\cos(M_y\omega_y y) \\
      \sin(\omega_x x)\sin(\omega_y y), & \sin(\omega_x x)\sin(2\omega_y y), & \cdots & \sin(M_x\omega_x x)\sin(M_y\omega_y y)
    \end{array}
    \right],
\end{equation}
where $\omega_x=2\pi/(x_r-x_l)$, $\omega_y=2\pi/(y_r-y_l)$, and $M_x,M_y$ are positive integers. The neural network then takes $(t,v(x,y))$ as input instead of $(t,x,y)$.

Since all components of the embedding~\eqref{eq:periodic_embedding_2d} are periodic in $x$ with period $x_r-x_l$ and periodic in $y$ with period $y_r-y_l$, the resulting network output is periodic in both spatial directions by construction. Moreover, if the activation functions are sufficiently smooth, the spatial derivatives obtained by automatic differentiation inherit the same periodicity from the sine and cosine features through the chain rule. Consequently, for the derivative orders required by the governing equation, including the zeroth-order case, we have
\begin{equation*}
  \partial_x^k \hat{u}(t,x_l,y;\theta) = \partial_x^k \hat{u}(t,x_r,y;\theta), \qquad \partial_y^k \hat{u}(t,x,y_l;\theta) = \partial_y^k \hat{u}(t,x,y_r;\theta),
\end{equation*}
where $k=0$ corresponds to the periodicity of the function value itself. Thus, the periodic boundary conditions are enforced exactly at both the function and derivative levels.

\section{Causal PINN}
\label{app:causal_pinn}

Conventional continuous-time PINNs typically minimize a residual loss aggregated over all space--time collocation points, which
implicitly treats residuals at all time instants as equally learnable throughout optimization. Wang et al.~\cite{wang2024Respecting} show that, for evolution equations, this practice can lead to a violation of temporal causality during training: residuals at later times may be reduced before the solution at earlier times is accurately resolved, allowing errors to propagate forward in time and trapping optimization in erroneous solutions. To enforce temporal precedence, they propose a simple reweighting of the residual loss that activates later-time residual terms only after earlier-time losses become sufficiently small.

Let $0=t_0<t_1<\cdots<t_{N_t}=T$ be a sequence of temporal points and $\{\mathbf{x}_j\}_{j=1}^{N_x}$ be spatial collocation points. Based on the residual-loss definition in~\eqref{eq:loss_terms} and the residual in~\eqref{eq:residual}, define the temporal residual loss at time $t_i$ as
\begin{equation}\label{eq:temporal_residual_loss}
  \mathcal{L}_{\mathrm{res}}(t_i,\theta)=\frac{1}{N_x}\sum_{j=1}^{N_x}\Big|\hat{u}_t(t_i,\mathbf{x}_j;\theta) +\mathcal{N}[\hat{u}](t_i,\mathbf{x}_j;\theta)\Big|^2.
\end{equation}
Using~\eqref{eq:temporal_residual_loss}, we incorporate the initial condition into the temporal ordering by treating the initial-condition loss as the first temporal constraint:
\begin{equation}\label{eq:temporal_loss}
  \mathcal{L}(t_0,\theta)=\omega_{\mathrm{ic}}\mathcal{L}_{\mathrm{ic}}(\theta), \qquad \mathcal{L}(t_i,\theta)=\omega_{\mathrm{res}}\,\mathcal{L}_{\mathrm{res}}(t_i,\theta)\quad (i\geqslant 1),
\end{equation}
where $\omega_{\mathrm{ic}},\omega_{\mathrm{res}}>0$ are the usual loss weights. In the numerical experiments considered in this work, boundary conditions are imposed as hard constraints. If boundary conditions are enforced weakly, one may incorporate them analogously by adding a term $\omega_{\mathrm{bc}}\mathcal{L}_{\mathrm{bc}}(t_i;\theta)$ to $\mathcal{L}(t_i;\theta)$ for the relevant time levels.

Based on the temporal loss~\eqref{eq:temporal_loss}, we define the weights $\{w_i\}_{i=0}^{N_t}$ by
\begin{equation}\label{eq:temporal_weights}
  w_0 = 1,\qquad w_i =\exp\Bigl(-\epsilon \sum_{j=0}^{i-1}\mathcal{L}(t_j,\theta)\Bigr),\quad i=1,2,\dots,N_t,
\end{equation}
where $\epsilon>0$ is a causality parameter controlling the steepness of the temporal gating. The resulting causal training objective is
\begin{equation}\label{eq:causal_loss}
  \mathcal{L}_{\mathrm{causal}}(\theta) = \frac{1}{N_t+1}\sum_{i=0}^{N_t} w_i\,\mathcal{L}(t_i,\theta).
\end{equation}
Because the weights \eqref{eq:temporal_weights} decay exponentially with the cumulative loss at earlier times, residual terms at $t_i$ receive negligible weight until all preceding constraints $\{\mathcal{L}(t_j,\theta)\}_{j<i}$ are sufficiently small; thus, optimization proceeds in a causality-respecting manner from early to late times. The PINN that is trained by minimizing \eqref{eq:causal_loss} is referred to as Causal PINN.

In the original implementation of Causal PINN, an annealing schedule for the causality parameter $\epsilon$ is used to gradually strengthen the temporal weighting. In our preliminary experiments, the annealing schedule did not consistently yield the best accuracy for the benchmark problems considered in this work. Therefore, to provide a stronger and more transparent baseline, we perform an independent sweep over fixed values of $\epsilon$ in the set
\begin{equation*}
  \epsilon \in \{10^{-4},10^{-3},10^{-2},10^{-1},10^{0},10^{1},10^{2},10^{3},10^{4}\},
\end{equation*}
while keeping all other hyperparameters unchanged. This procedure can be viewed as a one-dimensional parameter sweep over the causality parameter.

\bibliography{ref}

\end{document}